\documentclass[11pt, reqno]{amsart}
\usepackage{amssymb}
\usepackage{amscd}
\usepackage{verbatim,ifthen}
\usepackage{color}
\usepackage{latexsym}
\usepackage{mathrsfs}
\usepackage{wrapfig}
\usepackage{color}
\usepackage{bbm}
\usepackage{pgfplots}
\usepackage{multicol}
\usepackage{array}

\usepackage{amsthm} 

\newtheoremstyle{rightnum} 
  {3pt}                    
  {3pt}                    
  {\itshape}               
  {}                       
  {\bfseries}              
  {.}                      
  {.5em}                   
  {\thmname{#1}\thmnote{ \hfill(#3)}} 

\theoremstyle{rightnum}

\newcolumntype{M}[1]{>{\centering\arraybackslash}m{#1}}

\usepackage[colorlinks, linkcolor=black, citecolor=blue, linktocpage,backref=page]{hyperref}
\renewcommand*{\backref}[1]{}
\renewcommand*{\backrefalt}[4]{[{\tiny%
    \ifcase #1 Not cited.%
          \or Cited on page~#2.%
          \else Cited on pages #2.%
    \fi%
    }]}
\addtocontents{toc}{\setcounter{tocdepth}{1}}

\usepackage{amsmath}

\numberwithin{equation}{section}

\usepackage{eucal}
\usepackage{calc}  
\usepackage{enumitem} 
\usepackage{graphicx,wrapfig,lipsum}
\usepackage{etoolbox}
\usepackage{marginnote}
\usepackage{lipsum}
\makeatletter
\patchcmd{\@mn@margintest}{\@tempswafalse}{\@tempswatrue}{}{}
\patchcmd{\@mn@margintest}{\@tempswafalse}{\@tempswatrue}{}{}
\reversemarginpar 
\makeatother
\usepackage{scrextend}

\makeatletter
\DeclareRobustCommand\widecheck[1]{{\mathpalette\@widecheck{#1}}}
\def\@widecheck#1#2{%
    \setbox\z@\hbox{\m@th$#1#2$}%
    \setbox\tw@\hbox{\m@th$#1%
       \widehat{%
          \vrule\@width\z@\@height\ht\z@
          \vrule\@height\z@\@width\wd\z@}$}%
    \dp\tw@-\ht\z@
    \@tempdima\ht\z@ \advance\@tempdima2\ht\tw@ \divide\@tempdima\thr@@
    \setbox\tw@\hbox{%
       \raise\@tempdima\hbox{\scalebox{1}[-1]{\lower\@tempdima\box
\tw@}}}%
    {\ooalign{\box\tw@ \cr \box\z@}}}
\makeatother
\title{Hermitian metrics with vanishing second Chern Ricci curvature}
\author{Kyle Broder}
\address{The University of Queensland,  St. Lucia,  QLD 4067, Australia}
\email{k.broder@uq.edu.au}
\thanks{The first-named author is funded as a postdoctoral research fellow by an Australian Research Council Discovery Project Grant (DP220102530) held by the second-named author.}
\author{Artem Pulemotov}
\address{The University of Queensland,  St. Lucia,  QLD 4067, Australia}
\email{a.pulemotov@uq.edu.au}

\makeatother
\usepackage{xcolor}

\usepackage{soul}

\usepackage{calc}  
\usepackage{enumitem} 
\usepackage{tensor}
\usepackage{graphicx,wrapfig,lipsum}
\usepackage{etoolbox}
\usepackage{marginnote}
\usepackage{lipsum}
\makeatletter
\patchcmd{\@mn@margintest}{\@tempswafalse}{\@tempswatrue}{}{}
\patchcmd{\@mn@margintest}{\@tempswafalse}{\@tempswatrue}{}{}
\reversemarginpar 
\makeatother

\usepackage{array}
\newcolumntype{M}[1]{>{\centering\arraybackslash}m{#1}}

\usepackage{letltxmacro}
\makeatletter
\AtBeginDocument{%
  \@ifdefinable{\myorg@nameref}{%
    \LetLtxMacro\myorg@nameref\nameref
    \DeclareRobustCommand*{\nameref}[1]{%
      \emph{\myorg@nameref{#1}}%
    }%
  }%
}
\makeatother

\begin{document}

\maketitle 

\begin{abstract}
We describe a rigidity phenomenon exhibited by the second Chern
Ricci curvature of a Hermitian metric on a compact complex manifold. This yields a characterisation of second Chern Ricci-flat Hermitian metrics on several types of manifolds as well as a range of non-existence results for such metrics.

\end{abstract}

\section{Introduction and the main results}

\noindent The theory of Ricci-flat K\"ahler metrics has been a major component of complex geometry for almost half a century, pioneered by Yau's resolution of the Calabi conjecture \cite{Yau1976}. The past few decades have seen rapidly growing interest in the study of Hermitian generalizations of this theory on complex non-K\"ahler manifolds, with motivation coming from both mathematics and theoretical physics. Due to the presence of torsion in the Chern connection of a non-K\"ahler metric, the notions of Ricci curvature and, consequently, Ricci-flatness become ambiguous. In particular, there are four distinct \textit{Chern Ricci curvatures} for a Hermitian metric~$\omega$.

The first Chern Ricci curvature, denoted $\text{Ric}_{\omega}^{(1)}$, is given by the trace over the third and fourth indices of the Chern curvature tensor and is locally equal to $-\sqrt{-1} \partial \bar{\partial} \log(\omega^n)$. It represents the first (Bott--Chern) class (of the anti-canonical bundle) and is governed by a complex Monge--Amp\`ere equation (see, e.g., \cite{TosattiWeinkove, TosattiNonKahlerCY}). Compact non-K\"ahler Hermitian manifolds with $\text{Ric}_{\omega}^{(1)} =0$, called \textit{non-K\"ahler Calabi--Yau manifolds}, were systematically studied by Tosatti~\cite{TosattiNonKahlerCY} (see also \cite{AngellaCalamaiSpotti}). 

Our understanding of the second Chern Ricci curvature $\text{Ric}_{\omega}^{(2)}$, given by the trace over the first and second indices of the Chern curvature tensor,  is much less complete. This tensor does not have the cohomological nature that the first Chern Ricci curvature possesses. Its local expression is formidable and is not directly related to a complex Monge--Amp\`ere equation.  On the other hand,  the second Chern Ricci curvature plays essential roles in the Bochner technique~\cite{LiuYangGeometry, Lu, Royden, BroderSBC,BroderSBC2} and the Streets--Tian Hermitian curvature flow~\cite{StreetsTian}. Several remarkable properties of this tensor were established in ~\cite{LiuYangGeometry, LeeLi, YangSecondChernRicci}.   

Gauduchon--Ivanov~\cite{GauduchonIvanov} studied second Chern Einstein metrics (i.e., those $\omega$ for which $\text{Ric}_{\omega}^{(2)} = \lambda \omega$ with $\lambda \in \mathbf{R}$) on compact complex surfaces. They showed that all such metrics are K\"ahler--Einstein with the exception of the standard metric on the primary Hopf surface $\mathbf{S}^1 \times \mathbf{S}^3$, which has $\lambda >0$. As a consequence, second Chern Ricci-flat metrics on compact complex surfaces are Ricci-flat K\"ahler, and thus only exist on complex tori $\mathbf{C}^2/ \Lambda$, K3, and Enriques surfaces.

Angella--Calamai--Spotti~\cite{AngellaCalamaiSpotti} conducted a systematic investigation of second Chern Einstein metrics in all dimensions. They showed, among other things, that the possible signs of the Einstein constant correspond to the canonical bundle $K_X$ or the anti-canonical bundle $-K_X$ being pseudo-effective or unitary flat. Further,  second Chern Ricci-flat metrics have nonpositive Kodaira dimension $\text{kod}(X) \leq 0$.  

The remaining third and fourth Chern Ricci curvatures, $\text{Ric}_{\omega}^{(3)}$ and $\text{Ric}_{\omega}^{(4)}$, are given by the two remaining traces of the Chern curvature tensor. They are conjugate to each other, and they do not define real $(1,1)$--forms. In the framework developed by the first-named author and Tang~\cite{BroderTangAltered}, such curvatures provide altered variants. For a detailed explanation of this terminology, see~\cite[Remark~3.45]{BroderStanfield}.

In the present paper, we describe a rigidity phenomenon exhibited by the second Chern Ricci curvature and, most notably, by second Chern Ricci-flat metrics. This leads to several interesting corollaries, such as a complete understanding of the moduli space of second Chern Ricci-flat metrics on the torus. We also obtain several non-existence results and obstructions.  To state our first main theorem,  we remind the reader that the \textit{real bisectional curvature}~\cite{YangZhengRBC, LeeStreets} of a Hermitian metric $\omega$ is given by \begin{eqnarray*}
\text{RBC}_{\omega}(\xi) &:  = & \frac{1}{| \xi |_{\omega}^2} R_{\alpha \overline{\beta} \gamma \overline{\delta}} \xi^{\alpha \overline{\beta}} \xi^{\gamma \overline{\delta}},
\end{eqnarray*} where $R$ is the Chern curvature tensor and $\xi$ is a nonnegative Hermitian $(1,1)$--tensor.  This notion has played an essential role in the Schwarz lemma for holomorphic maps between Hermitian manifolds ~\cite{YangZhengRBC, LeeStreets,BroderSBC,BroderSBC2}.

\subsection*{Theorem~1.1}\label{Theorem_Neg}
Let $(X,\widetilde{\omega})$ be a compact Hermitian manifold with $\text{RBC}_{\widetilde{\omega}} \leq 0$.  If the equality $\text{Ric}_{\omega}^{(2)}=0$ holds for some Hermitian metric $\omega$ on~$X$, then $\omega$ has the same Chern connection as~$\widetilde\omega$. If there is a point where $\text{RBC}_{\widetilde{\omega}}<0$, then there are no metrics on $X$ with vanishing second Chern Ricci curvature. \\

For K\"ahler metrics, the real bisectional curvature is comparable to (in the sense that it always has the same sign as) the more familiar \textit{holomorphic sectional curvature} $$\text{HSC}_{\omega}(v) \ : = \ \frac{1}{| v |_{\omega}^4} R_{\alpha \overline{\beta} \gamma \overline{\delta}} v^{\alpha} \overline{v}^{\beta} v^{\gamma} \overline{v}^{\delta},$$ where $v \in T^{1,0}X$. This holds more generally, for the class of K\"ahler-like metrics~\cite{YangZhengCurvature},  described by the Chern curvature tensor having the same symmetries as that of a K\"ahler metric.  Since the K\"ahler condition $d \omega =0$ is equivalent to the vanishing of the Chern torsion, the property of being K\"ahler depends only on the Chern connection. The same holds for the property of being K\"ahler-like. With these observations at hand, we can derive the following from \nameref{Theorem_Neg}.

\subsection*{Corollary~1.2}\label{corollary_met_type}
Assume that a compact complex manifold $X$ admits a K\"ahler (respectively,  K\"ahler-like) metric with $\text{HSC}_{\widetilde{\omega}} \leq 0$.  Then every Hermitian metric $\omega$ on $X$ satisfying $\text{Ric}_{\omega}^{(2)}=0$ is K\"ahler (respectively, K\"ahler-like). If such a metric exists,  the manifold $X$ is Calabi--Yau in the sense that the canonical bundle is holomorphically torsion. \\

Berger~\cite{BergerHBC} showed that any Ricci-flat Riemannian metric on the torus $\mathbf{T}^4$ is flat. This was extended to all dimensions by Fischer--Wolf~\cite{FischerWolf}.  Since complex tori $\mathbf{T}^n$ admit flat K\"ahler metrics, \nameref{corollary_met_type}  yields a complete description of the moduli space of Hermitian metrics with vanishing second Chern Ricci curvature on~$\mathbf T^n$. 

\subsection*{Corollary~1.3}\label{corollary_torus}
Every second Chern Ricci-flat metric on a complex torus is flat.\\

The Alekseevsky--Kimelfeld theorem~\cite[Theorem~7.61]{AB87} asserts that every homogeneous Ricci-flat Riemannian metric must be flat. It is noteworthy that the result of Fischer--Wolf and \nameref{corollary_torus}, while similar to this theorem in spirit, hold without any symmetry assumptions.

K\"ahler metrics being the most prominent example, several of the types of Hermitian metrics that arise in complex geometry can be described exclusively in terms of the Chern connection. For instance, this is the case for balanced metrics. Thus, from \nameref{Theorem_Neg}, we obtain the following analogue of \nameref{corollary_met_type}.

\subsection*{Corollary~1.4}\label{corollary_met_type2}
Assume that a compact complex manifold $X$ admits a balanced metric with $\text{RBC}_{\widetilde{\omega}} \leq 0$.  Then every Hermitian metric $\omega$ on $X$ satisfying $\text{Ric}_{\omega}^{(2)}=0$ is balanced.  \\

We will not endeavor to provide an exhaustive list of the types of metrics determined by the Chern connection that have appeared in the literature. It is obvious, however, that \nameref{corollary_met_type2} extends to all of them.   

Our second main theorem is a counterpart to \nameref{Theorem_Neg} for nonnegatively curved metrics. To state it, we remind the reader of the \emph{Schwarz bisectional curvature} that was introduced by the first-named author ~\cite{BroderSBC, BroderSBC2}: 
$$\text{SBC}_{\omega}(\xi) \ : = \  R_{\alpha \overline{\beta} \gamma \overline{\delta}} \xi^{\alpha \overline{\beta}} (\xi^{-1})^{\gamma \overline{\delta}},$$ where $\xi \in \Lambda_X^{1,1}$ is a positive-definite Hermitian $(1,1)$--tensor with inverse $\xi^{-1}$.
The Schwarz bisectional curvature provides an analogue of the real bisectional curvature for the Aubin--Yau formulation of the Schwarz lemma.

\subsection*{Theorem~1.5}\label{Theorem_Pos}
Let $(X,\widetilde{\omega})$ be a compact Hermitian manifold with $\text{SBC}_{\widetilde{\omega}} \geq 0$. If $\text{Ric}_{\omega}^{(2)}=0$ for some Hermitian metric $\omega$ on~$X$, then $\omega$ has the same Chern connection as~$\widetilde\omega$. If there is a point where $\text{SBC}_{\widetilde{\omega}}>0$, then there are no metrics on $X$ with vanishing second Chern Ricci curvature.\\

The strength of Schwarz bisectional curvature remains unknown, but it is dominated by the holomorphic bisectional curvature.  As a consequence, we obtain the following.

\subsection*{Corollary~1.6}\label{Cor17}
There are no metrics with vanishing second Chern Ricci curvature on compact Hermitian symmetric spaces. \\

The proofs of \nameref{Theorem_Neg} and \nameref{Theorem_Pos} rely on the Chern--Lu ~\cite{Lu,Royden,YangZhengRBC} and Aubin--Yau ~\cite{Yau1976,BroderSBC,BroderSBC2} incarnations of the Schwarz lemma, respectively.  The way we employ these results is inspired by the work of DeTurck--Koiso~\cite{DeTurckKoiso} concerning the prescribed Ricci curvature problem on Riemannian manifolds. In fact, as a byproduct of our techniques, we obtain a uniqueness theorem for Hermitian metrics with prescribed second Chern Ricci curvature; see Section~\ref{DTKSection}.  Interestingly, thanks to the underlying Hermitian structure, this theorem requires weaker assumptions than its (real) Riemannian counterpart.

\subsection*{Acknowledgements}
The authors would like to thank Ben Andrews, Ramiro Lafuente, James Stanfield, and Wolfgang Ziller for valuable discussions. Part of this work was completed while the authors were visiting Ground Zero. We would like to thank them for their hospitality and express our gratitude to G\"urkan Baydogan for stimulating conversations that helped us grapple with some of the ideas in this paper. 

\section{Preliminaries}
\noindent Throughout we will denote by $X$ a complex manifold of (complex) dimension~$n$. Let $\mathcal J$ be the underlying complex structure. It splits the complexified tangent bundle $T^{\mathbf{C}}X : = TX \otimes_{\mathbf{R}} \mathbf{C}$ into a sum of eigenbundles $T^{1,0}X \oplus T^{0,1}X$, where $T^{1,0}X$ is the $\sqrt{-1}$--eigenbundle of $(1,0)$--tangent vectors $u - \sqrt{-1} \mathcal{J} u$ and $T^{0,1}X$ is the $-\sqrt{-1}$--eigenbundle of $(0,1)$--tangent vectors $u + \sqrt{-1} \mathcal{J} u$.  The splitting of $T^{\mathbf{C}} X$ induces a splitting of the bundle of $k$--forms: $\Lambda_X^k \ \simeq \ \bigoplus_{p+q=k} \Lambda_X^{p,q}$. We denote by $K_X : = \Lambda_X^{n,0}$  the canonical bundle, whose dual is the anti-canonical bundle~$-K_X$. The exterior derivative splits as $d = \partial + \bar{\partial}$ with $\partial : \Omega_X^{p,q} \to \Omega_X^{p+1, q}$ and $\bar{\partial} : \Omega_X^{p,q} \to \Omega_X^{p,q+1}$. 

A Riemannian metric $g$ on a complex manifold $(X,\mathcal{J})$ is said to be \textit{Hermitian} if $$g(\mathcal{J}\cdot, \mathcal{J} \cdot) \ = \ g(\cdot, \cdot).$$ We write $\omega=\omega_g(\cdot, \cdot) : = g(\mathcal{J} \cdot, \cdot)$ for the associated $2$--form and follow the general practice of referring to $\omega$ as a Hermitian metric, often suppressing any mention of $g$ or~$\mathcal{J}$. Several distinguished metric types are given in terms of the associated $2$--form: a Hermitian metric $\omega$ is \textit{K\"ahler} if $d\omega=0$, and  \textit{balanced} if $d \omega^{n-1}=0$. These structures can be encoded in the properties of the \textit{Chern connection}---the unique Hermitian connection such that $\nabla^{0,1} = \bar{\partial}$.  Indeed, the K\"ahler metrics are characterized by the vanishing of the Chern torsion, while the balanced metrics are those for which the torsion $(1,0)$--form $\tau = \sum_{i,k} T_{ik}^k e^i$ vanishes, where $e^i$ is any local coframe. 

Let $R$ denote the Chern curvature tensor. In a local frame $e_k$ for $T^{1,0}X$, the only non-zero components of $R$ are $R(e_i, \overline{e}_j) e_k = - \partial_{\overline{j}} \Gamma_{ik}^p e_p =: R_{i \overline{j} k}{}^{p} e_{p}$. We use the Hermitian metric to lower the index, writing $R_{i \overline{j} k \overline{\ell}} := R_{i \overline{j} k}{}^p g_{p \overline{\ell}}$. Further, we note that the Chern curvature tensor has the conjugate symmetry $R_{\overline{j} i \overline{\ell} k} = \overline{R_{i \overline{j} k \overline{\ell}}}$.  For a general Hermitian non-K\"ahler metric,  the Chern curvature tensor will not possess the symmetries of the Riemannian curvature tensor.  

While the sectional curvature of the Chern connection is meaningful, the more natural counterpart of the sectional curvature in complex geometry is the \textit{holomorphic bisectional curvature} \begin{eqnarray}\label{HBC}
\text{HBC}_{\omega}(u,v) & : = & \frac{1}{| u |_{\omega}^2 | v |_{\omega}^2} R_{\alpha \overline{\beta} \gamma \overline{\delta}} u^{\alpha} \overline{u}^{\beta} v^{\gamma} \overline{v}^{\delta},
\end{eqnarray} where $u,v \in T^{1,0}X$ are $(1,0)$--tangent vectors. Compact K\"ahler manifolds with $\text{HBC}_{\omega} \geq 0$ are Hermitian symmetric spaces of compact type by Mok's resolution of the generalized Frankel conjecture~\cite{MokFrankel}.

\section{Nonpositively curved metrics}

\noindent In this section, we proof \nameref{Theorem_Neg} and \nameref{corollary_met_type}. We also discuss several additional corollaries of these results. The main technique we use is the Schwarz lemma, which is the product of combining the maximum principle with the Bochner formula applied to $| \partial f |^2$ for a holomorphic map~$f$.  We, therefore, start by recalling the main Schwarz lemma calculation that goes back to Lu~\cite{Lu}. Given a holomorphic map  $f : (X, \omega) \to (X, \widetilde{\omega})$ and local holomorphic coordinates $(z_1,..., z_n)$ centered at a point $p \in X$, let us write $f_i^{\alpha} : = \frac{\partial f^{\alpha}}{\partial z_i}$,  where $f = (f^1, ..., f^n)$. The complex Laplacian $\Delta_{\omega}$ with respect to the metric $\omega$ is defined by $\Delta_{\omega} : = \text{tr}_{\omega} \sqrt{-1} \partial \bar{\partial}$.

\begin{proof}[Proof of Theorem~1.1]
Let $f : (X, \omega) \to (X, \widetilde{\omega})$ be a holomorphic map. We will eventually take $f = \text{id}$,  but it is more transparent to work with a general holomorphic map,  restricting to the identity map at the end of the proof. Let $g$ and $\widetilde{g}$ be the metrics underlying $\omega$ and $\widetilde{\omega}$, respectively. Denote by $\nabla$ and $\widetilde{\nabla}$ the Chern connections of $\omega$ and $\widetilde{\omega}$. Write $\Gamma_{ij}^k = g^{k \overline{\ell}} \partial_i g_{j \overline{\ell}}$ and $\widetilde{\Gamma}_{\alpha \beta}^{\gamma} = \widetilde{g}^{\gamma \overline{\delta}} \partial_{\alpha} \widetilde{g}_{\beta \overline{\delta}}$ for the corresponding Christoffel symbols.  The differential $\partial f$ is a holomorphic section of $\Omega_X^{1,0} \otimes f^{\ast} T^{1,0}X$.  Hence, if we write $\widehat{\nabla}$ for the connection on $\Omega_X^{1,0} \otimes f^{\ast} T^{1,0}X$ induced from $\nabla$ and $\widetilde{\nabla}$, then in any local frame,  \begin{eqnarray}\label{DiffChristoffel}
(\widehat{\nabla}_k \partial f)_{\ell}^{\alpha} &=& f_{k\ell}^{\alpha} + \widetilde{\Gamma}_{\gamma \rho}^{\alpha} f_k^{\gamma} f_{\ell}^{\rho} - \Gamma_{k\ell}^m f_m^{\alpha}.
\end{eqnarray} A calculation going back to Lu ~\cite{Lu} shows that \begin{eqnarray}\label{CoordinateSchwarzFormula}
\Delta_{\omega} | \partial f |^2  &=& | \widehat{\nabla} \partial f |^2 + \text{Ric}_{k \overline{\ell}}^{(2)} g^{k \overline{q}} g^{p \overline{\ell}} \widetilde{g}_{\alpha \overline{\beta}} f_p^{\alpha} \overline{f_q^{\beta}} - g^{i \overline{j}} g^{p \overline{q}} \widetilde{R}_{\alpha \overline{\beta} \gamma \overline{\delta}} f_i^{\alpha} \overline{f_j^{\beta}}f_p^{\gamma} \overline{f_q^{\delta}},
\end{eqnarray} where $\text{Ric}^{(2)}_{k \overline{\ell}}$ are the components of the second Chern Ricci curvature of $\omega$ and $\widetilde{R}$ is the Chern curvature tensor of $\widetilde{\omega}$.  Since $\text{Ric}^{(2)}_{\omega}=0$ and $\text{RBC}_{\widetilde{\omega}} \leq   0$,  it follows that \begin{eqnarray}\label{Estimate1}
\Delta_{\omega} | \partial f |^2 &=& | \widehat{\nabla} \partial f |^2 - g^{i \overline{j}} g^{p \overline{q}} \widetilde{R}_{\alpha \overline{\beta} \gamma \overline{\delta}} f_i^{\alpha} \overline{f_j^{\beta}}f_p^{\gamma} \overline{f_q^{\delta}}  \ \geq \ | \widehat{\nabla} \partial f |^2.
\end{eqnarray} Because $X$ is compact, \begin{eqnarray}\label{IntDiv}
0 \ = \ \int_X \Delta_{\omega} | \partial f |^2 \omega^n & \geq & \int_X | \widehat{\nabla} \partial f |^2 \omega^n, 
\end{eqnarray} and therefore, $| \widehat{\nabla} \partial f |^2 =0$. Taking $f = \text{id}$ in \eqref{DiffChristoffel}, we see that  \begin{eqnarray*}
(\widehat{\nabla}_k \partial f)_{\ell}^m &=& \widetilde{\Gamma}_{\gamma \rho}^{\alpha} \delta_k^{\gamma} \delta_{\ell}^{\rho} - \Gamma_{k\ell}^m \delta_m^{\alpha} \ = \ \widetilde{\Gamma}_{k \ell}^m - \Gamma_{k\ell}^m,
\end{eqnarray*} and so the Christoffel symbols of $\omega$ and $\widetilde{\omega}$ coincide.  It follows from \eqref{Estimate1} and \eqref{IntDiv} that if $\text{RBC}_{\widetilde{\omega}} < 0$ at some point,  then we arrive at a contradiction, showing that there are no second Chern Ricci-flat metrics on compact manifolds with quasi-negative real bisectional curvature.
\end{proof}

\begin{proof}[Proof of Corollary~1.2]
Let $f : (X, \omega) \to (X, \widetilde{\omega})$ be a holomorphic map.  Suppose that $\widetilde{\omega}$ is K\"ahler-like.  By Royden's polarisation argument~\cite[p. 552]{Royden}, $\text{HSC}_{\widetilde{\omega}} \leq 0$ (respectively, $\text{HSC}_{\widetilde{\omega}} <0$) implies that $\text{RBC}_{\widetilde{\omega}} \leq 0$ (respectively, $\text{RBC}_{\widetilde{\omega}} <0$). Hence, we can apply \nameref{Theorem_Neg} to obtain the equality of the Chern connections in the $\text{HSC}_{\widetilde{\omega}}\leq 0$ case, and non-existence if $\text{HSC}_{\widetilde{\omega}} <0$.
\end{proof}

\subsection*{Remark~3.1}\label{Remark33}
Instead of assuming that $\text{RBC}_{\widetilde\omega}\le0$ and $\text{Ric}_\omega^{(2)}=0$ in \nameref{Theorem_Neg}, one may assume that $\text{RBC}_{\widetilde{\omega}} \leq c$ and $\text{Ric}_{\omega}^{(2)} = c \widetilde{\omega}$ for some constant $c\in\mathbf R$. Only minor changes to the proof are required. \nameref{corollary_met_type} admits a similar generalisation.\\

In the case where the real bisectional curvature vanishes identically, Yang--Zheng ~\cite{YangZhengRBC} showed that the metric must be balanced. A metric that is simultaneously pluriclosed and balanced must be K\"ahler. This is clear from tracing the formula for the second Chern Ricci curvature (see, e.g., \cite[Theorem 3.4]{BroderStanfield}).

\subsection*{Corollary~3.2}
Let $(X, \widetilde{\omega})$ be a compact Hermitian manifold with real bisectional curvature vanishing identically. Every second Chern Ricci-flat metric on $X$ is balanced.  If $X$ is non-K\"ahler, there are no second Chern Ricci-flat pluriclosed metrics on~$X$. \\

It is clear from \cite[Theorem 2.9]{BroderTangAltered} that the assertion of this corollary also holds for compact Hermitian manifolds whose altered real bisectional curvature (in the sense of~\cite{BroderTangAltered}) vanishes identically.

\section{Nonnegatively curved metrics}

\noindent In this section, we prove \nameref{Theorem_Pos} and \nameref{Cor17}. The main technique is the Hermitian Aubin--Yau inequality established by the first-named author in~\cite{BroderSBC} (building from~\cite{Aubin,Yau1976}). 

\begin{proof}[Proof of Theorem~1.5.]
As in the proofs of \nameref{Theorem_Neg} and \nameref{corollary_met_type}, we will work with a general biholomorphic map $f : (X, \omega) \to (X, \widetilde{\omega})$ and restrict to the identity map at the end. We also maintain the notation used in the proof of \nameref{Theorem_Neg}.  From~\cite{BroderSBC},  we have
\begin{align*}
\Delta_{\widetilde\omega} | \partial f |^2 \circ f^{-1} = | \widehat{\nabla} \partial f \circ f^{-1} |^2 &- \widetilde{\text{Ric}}_{\alpha \overline{\beta}}^{(2)} g^{i \overline{j}} f_i^{\alpha} \overline{f_j^{\beta}} 
\\
&+ R_{k \overline{\ell} p \overline{q}} g^{i \overline{q}} g^{p \overline{j}} \widetilde g^{\gamma \overline{\delta}} \widetilde g_{\alpha \overline{\beta}} f_i^{\alpha} \overline{f_j^{\beta}} (f^{-1})_{\gamma}^k \overline{(f^{-1})_{\delta}^{\ell}}.
\end{align*} The source curvature term $$R_{k \overline{\ell} p \overline{q}} g^{\gamma \overline{\delta}} g_{\alpha \overline{\beta}} \widetilde g^{i \overline{q}} \widetilde g^{p \overline{j}} f_i^{\alpha} \overline{f_j^{\beta}} (f^{-1})_{\gamma}^k \overline{(f^{-1})_{\delta}^{\ell}}$$
is controlled precisely by the Schwarz bisectional curvature,  and hence
\begin{eqnarray*}
\Delta_{\widetilde\omega} | \partial f |^2 \circ f^{-1} & \geq & | \widehat{\nabla} \partial f \circ f^{-1} |^2.
\end{eqnarray*}
Integrating this inequality, we see that $| \widehat{\nabla} \partial f \circ f^{-1} |^2 = 0$.  Taking $f$ to be the identity map shows that the Chern connections coincide. Moreover, if $\text{SBC}_{\widetilde{\omega}} >0$ at one point, we get a contradiction in the same manner as in the proof of \nameref{Theorem_Neg}.
\end{proof}

The strength of the Schwarz bisectional curvature remains largely unknown. However, it is clear that the holomorphic bisectional curvature dominates the Schwarz bisectional curvature in the sense that the sign of the former determines the same sign on the latter. This enables us to prove \nameref{Cor17}.

\begin{proof}[Proof of Corollary~1.6]
A compact Hermitian symmetric space $X$ admits a K\"ahler metric $\widetilde\omega$ with $\text{HBC}_{\widetilde{\omega}} \geq 0$ (see, e.g., \cite{MokZhong}).  By \nameref{Theorem_Pos},  since the holomorphic bisectional curvature dominates the Schwarz bisectional curvature, every second Chern Ricci-flat metric on $X$ must be K\"ahler.  If such a metric existed,  the canonical bundle would be holomorphically torsion, which is impossible.
\end{proof}

Let us remark that a second Chern Ricci-flat metric was found on the (non-compact) Snow manifold $S5$ in \cite[$\S 3.3.4$]{AngellaCalamaiSpotti}.

\section{Hermitian metrics with prescribed second Chern Ricci curvature}\label{DTKSection}

\noindent In the K\"ahler category, the prescribed Ricci curvature problem is settled by Yau's resolution of the Calabi conjecture \cite{Yau1976}. In the Riemannian category, however, this problem is largely unresolved. Numerous results have appeared in the literature over the past four decades (see the survey~\cite{BP19} and more recent works such as~\cite{BK20,LWa,PulemotovZiller}). One of the first uniqueness theorems for Riemannian metrics with prescribed Ricci curvature was produced by DeTurck--Koiso~\cite{DeTurckKoiso}. According to this theorem, if $(M, \bar{g})$ is a compact Einstein manifold with nonnegative sectional curvature and $\text{Ric}_{\bar{g}} = \bar{g}$, then every Riemannian metric $g$ on $M$ with $\text{Ric}_g = \text{Ric}_{\bar{g}}$ has the same Levi-Civita connection as $\bar{g}$. Analysis of K\"ahler metrics on the Wallach space $\text{SU}(3)/\mathbf{T}^2$ demonstrates that the result may fail without the nonnegativity assumption on the sectional curvature (see~\cite[Example~1]{PulemotovZiller}). A similar uniqueness theorem holds in the Hermitian category for the second Chern Ricci curvature. In this setting, we require only the nonnegativity of the holomorphic bisectional curvature. 

\subsection*{Theorem 5.1}\label{DeTurckKoisoTrick}
Let $(X, \widetilde{\omega})$ be a compact Hermitian manifold with $\text{Ric}_{\widetilde{\omega}}^{(2)} = \widetilde{\omega}$ and $\text{HBC}_{\widetilde{\omega}} \geq 0$. If $\omega$ is a Hermitian metric on $X$ with $\text{Ric}_{\omega}^{(2)} = \text{Ric}_{\widetilde{\omega}}^{(2)}$, then $\omega$ and $\widetilde{\omega}$ have the same Chern connection. 

\begin{proof}
The assumptions on $\widetilde{\omega}$ ensure that the real bisectional curvature of $\widetilde{\omega}$ is bounded above. To see this, denote by $\mathcal{H}_X$ the space of Hermitian $(0,2)$--tensors on $X$. Let $\widetilde{g}$ be the metric underlying $\widetilde{\omega}$.  The real bisectional curvature $\text{RBC}_{\widetilde{\omega}}$ is the quadratic form associated to the complex curvature operator $$\mathfrak{K} : \mathcal{H}_X \to \mathcal{H}_X, \hspace{1cm} \mathfrak{K}(\xi)_{k \overline{\ell}} \ : = \  \widetilde{R}_{i \overline{j} k \overline{\ell}} \widetilde{g}^{p \overline{j}} \widetilde{g}^{i \overline{q}} \xi_{p \overline{q}},$$ where $\widetilde{R}$ is the Chern curvature tensor of $\widetilde{\omega}$. Hence, an upper bound on $\text{RBC}_{\widetilde{\omega}}$ is given by an upper bound on the largest eigenvalue of $\mathfrak{K}$. Since $\widetilde{g}$ is second Chern Einstein, it is an eigenvector of $\mathfrak{K}$ with unit eigenvalue: \begin{eqnarray*}
    \mathfrak{K}(\widetilde{g}_{p \overline{q}}) &=& \widetilde{R}_{i \overline{j} k \overline{\ell}} \widetilde{g}^{p \overline{j}} \widetilde{g}^{i \overline{q}} \widetilde{g}_{p\bar{q}} \ = \ \widetilde{g}^{i \overline{j}} \widetilde{R}_{i \overline{j} k \overline{\ell}} \ = \ \widetilde{\text{Ric}}_{k \overline{\ell}}^{(2)} \ = \ \widetilde{g}_{k \overline{\ell}}.
\end{eqnarray*} Let $\zeta$ be an eigenvector of $\mathfrak{K}$, distinct from $\widetilde g$, with eigenvalue $\mu$. Since $\mathfrak{K}$ is Hermitian, the distinct eigenvectors of $\mathfrak{K}$ are orthogonal. Hence, we may assume that $\text{tr}_{\widetilde{g}}(\zeta) =0$ and $\lambda_1 = \max_k  \lambda_k  >0$. Fix a point $x_0\in X$ and choose coordinates such that $\widetilde{g}_{i \overline{j}} = \delta_{ij}$ and $\zeta_{i \overline{j}} = \lambda_i \delta_{ij}$ at~$x_0$. Computing at this point, we find
\begin{eqnarray*}
\mu \lambda_1 \ = \ \mu \zeta_{1 \overline{1}} \ = \ \mathfrak{K}(\zeta)_{1 \overline{1}} &=& \sum \widetilde{R}_{i \overline{j} 1 \overline{1}}  \widetilde{g}^{i\overline{q}} \widetilde{g}^{p \overline{j}} \zeta_{p \overline{q}} \ = \ \sum_p \widetilde{R}_{p \overline{p} 1 \overline{1}} \lambda_p.
\end{eqnarray*}
Hence, \begin{eqnarray*}
\mu  \ = \  \frac{1}{\lambda_1} \sum_p \widetilde{R}_{p \overline{p} 1 \overline{1}} \lambda_p &=& \frac{1}{\lambda_1} \sum_p (\widetilde{R}_{p \overline{p} 1 \overline{1}} - \mathcal{B}_{\min}) \lambda_p \\
& \leq & \frac{1}{\lambda_1} \sum_p (\widetilde{R}_{p \overline{p} 1 \overline{1}} - \mathcal{B}_{\min}) \lambda_1 \ = \ 1 - n \mathcal{B}_{\min},
\end{eqnarray*} where $\mathcal{B}_{\min}$ denotes the minimum of the holomorphic bisectional curvature. It follows that the eigenvalues of $\mathfrak{K}$ are bounded above by $\max\{ 1, 1 - n \mathcal{B}_{\min} \}$. Therefore, if $\mathcal{B}_{\min} \geq 0$, we must have $\text{RBC}_{\widetilde{\omega}} \leq 1$. 

Let now $f : (X, \omega) \to (X, \widetilde{\omega})$ be a holomorphic map.  Since $\text{Ric}_{\omega}^{(2)} = \text{Ric}_{\widetilde{\omega}}^{(2)} = \widetilde{\omega}$, we see that \begin{eqnarray*}
\text{Ric}_{k \overline{\ell}}^{(2)} g^{k \overline{q}} g^{p \overline{\ell}} \widetilde{g}_{\alpha \overline{\beta}} f_p^{\alpha} \overline{f_q^{\beta}} &=& f_k^{\gamma} \overline{f_{\ell}^{\delta}} \widetilde{g}_{\gamma \overline{\delta}} g^{k \overline{q}} g^{p \overline{\ell}} \widetilde{g}_{\alpha \overline{\beta}} f_p^{\alpha} \overline{f_q^{\beta}} \ = \ | \partial f |^4.
\end{eqnarray*} Hence, using the inequality $\text{RBC}_{\widetilde{\omega}} \leq 1$ and the Chern--Lu formula~\eqref{CoordinateSchwarzFormula}, we find
\begin{eqnarray*}
    \Delta_{\omega} | \partial f |^2 &=& | \widehat{\nabla} \partial f |^2 + | \partial f |^4  - g^{i \overline{j}} g^{p \overline{q}} \widetilde{R}_{\alpha \overline{\beta} \gamma \overline{\delta}} f_i^{\alpha} \overline{f_j^{\beta}}f_p^{\gamma} \overline{f_q^{\delta}} \ \geq \ | \widehat{\nabla} \partial f |^2.
\end{eqnarray*} Since $X$ is compact, by integrating the above inequality with $f = \text{id}$ as in the proof of \nameref{Theorem_Neg}, we obtain the desired conclusion.
\end{proof}

The main theorem of~\cite{MokZhong} states that a compact K\"ahler manifold $(X,\widetilde{\omega})$ with $\text{Ric}_{\widetilde{\omega}}>0$ and $\text{HBC}_{\widetilde{\omega}} \geq 0$ must be biholomorphically isometric to a compact Hermitian symmetric space. It is interesting to compare \nameref{DeTurckKoisoTrick} to this result. The conditions on the curvature in \nameref{DeTurckKoisoTrick} are considerably stronger: we require $\text{Ric}_\omega^{(2)}=\text{Ric}_{\widetilde\omega}^{(2)}$ instead of just positivity. However, the gain is that \nameref{DeTurckKoisoTrick} does not require the manifold $X$ to be K\"ahler. For instance, in (complex) dimension $2$,  there are two metrics $\widetilde{\omega}$ that satisfy the hypotheses: the (K\"ahler) Fubini--Study metric on $\mathbf{P}^2$ and the (non-K\"ahler) standard metric on a primary Hopf surface $\mathbf{S}^1 \times \mathbf{S}^3$ (see~\cite{GauduchonIvanov}).

\subsection*{Corollary~5.2}
Let $\mathbf{S}^1 \times \mathbf{S}^3$ be a primary Hopf surface endowed with its standard metric~$\omega_0$.  If $\omega$ is a Hermitian metric on $\mathbf{S}^1 \times \mathbf{S}^3$ with $\text{Ric}_{\omega}^{(2)} = \text{Ric}_{\omega_0}^{(2)}$, then $\omega$ has the same Chern connection as $\omega_0$.

\subsection*{Remark~5.3}\label{Remark_HermBetter}
It is interesting to compare our use of the Bochner technique to that of DeTurck--Koiso \cite{DeTurckKoiso}. The identity map $\text{id} : (M,g) \to (M,\widetilde{g})$ of Riemannian manifolds is harmonic only if $g^{jk}(\Gamma_{jk}^i - \widetilde{\Gamma}_{jk}^i)=0$. From the Bianchi identity,  this is satisfied by $\text{id} : (M, g) \to (M, \text{Ric}_g)$ if $\text{Ric}_{g}$ is positive definite, and this is the situation that DeTurck--Koiso consider. In the Hermitian category, the identity map is holomorphic, and hence, the Bochner technique can be applied without these cumbersome assumptions on the Christoffel symbols or the positivity of the Ricci curvature.

\subsection*{Remark~5.4}\label{Remark_Ric1_Ric2}
The first Chern Ricci curvature of a Hermitian metric is determined completely by the Chern connection of the metric. However, since $\text{Ric}^{(2)}_\omega= g^{i \overline{j}} R_{i \overline{j}k \overline{\ell}}$, it appears that two Hermitian metrics may have distinct second Chern Ricci curvatures even if their Chern connections coincide. For this reason, we do not claim that $\text{Ric}^{(2)}_\omega=\text{Ric}^{(2)}_{\widetilde\omega}$ in \nameref{Theorem_Neg} and \nameref{Theorem_Pos}. Similarly, it appears that the pluriclosed condition is not determined entirely by the Chern connection. Indeed, let $\Gamma_{ij}^k = g^{k \overline{\ell}} \partial_i g_{j \overline{\ell}}$ be the Christoffel symbols of the Chern connection. The pluriclosed condition $\sqrt{-1} \partial \bar{\partial} \omega=0$ is equivalent to \begin{eqnarray*}
    g_{p \overline{\ell}} \left( R_{k \bar{j} i}{}^p - R_{i \bar{j} k}{}^p \right) + g_{p \bar{j}} \left( R_{i \bar{\ell} k}{}^p - R_{k \bar{\ell} i}{}^p \right)   + 4 g_{p \bar{q}}  \Gamma_{ik}^p \overline{\Gamma_{j\ell}^q} &=& 0.
\end{eqnarray*} This cannot be expressed in terms of the Christoffel symbols alone.

\end{document}